\documentclass{amsart}
\usepackage[all]{xy}
\usepackage{amssymb,url}
\newcommand{\calB}{\mathcal{B}}
\newcommand{\calC}{\mathcal{C}}
\newcommand{\calG}{\mathcal{G}}
\newcommand{\rtG}{\operatorname{rt}\nolimits\calG}
\newcommand{\calI}{\mathcal{I}}
\renewcommand{\L}{\Lambda}
\newcommand{\extto}{\xrightarrow}
\newcommand{\tip}{\operatorname{tip}\nolimits}
\newcommand{\tippath}{\operatorname{tippath}\nolimits}
\newcommand{\tipcoord}{\operatorname{tipcoord}\nolimits}
\newcommand{\coefftip}{\operatorname{coefficient\ of\ tip}\nolimits}
\newcommand{\rrad}{\mathfrak{r}}
\newcommand{\Ker}{\operatorname{Ker}\nolimits}
\renewcommand{\Im}{\operatorname{Im}\nolimits}
\newcommand{\id}{\operatorname{id}\nolimits}
\newcommand{\norm}{\operatorname{\textbf{NormalForm}}\nolimits}

\newtheorem{lem}{Lemma}[section]
\newtheorem{prop}[lem]{Proposition}
\newtheorem{thm}[lem]{Theorem}
\newtheorem{cor}[lem]{Corollary}
\theoremstyle{definition}

\newtheorem{example}[lem]{Example}

\title{An algorithmic approach to resolutions}
\author[Green]{Edward L. Green$^\dagger$}
\thanks{$^\dagger$Partially supported by a grant from the NSA and 
travel support from the Research Council of Norway}
\address{Edward L. Green\\
Department of Mathematics\\ 
Virginia Tech\\
Blacksburg, VA 24061--0123\\
USA}
\email{green@math.vt.edu}
\author[Solberg]{\O yvind Solberg$^\ddagger$}
\thanks{$^\ddagger$Partially supported by the Research Council of Norway}
\address{\O yvind Solberg\\
Institutt for matematiske fag\\
NTNU\\ 
N--7491 Trondheim\\ 
Norway}
\email{oyvinso@math.ntnu.no}

\begin{document}
\begin{abstract}
We provide an algorithmic method for constructing projective
resolutions of modules over quotients of path algebras. This algorithm
is modified to construct minimal projective resolutions of linear
modules over Koszul algebras. 
\end{abstract}

\maketitle

\section{Introduction}
Projective resolutions play an important role in homological
algebra. The existence of algorithmic methods has led to programs that
construct projective resolutions of modules. For commutative rings we
mention CoCoa System, Faug\`ere's GB, Macaulay, and Singular
\cite{Co,Fa,Ma,Si}. On the other hand, for non-commutative rings there
are fewer choices. For group algebras there are programs written for
MAGMA by Jon Carlson \cite{MA}.  The program Bergman \cite{Be}
provides projective resolutions for quotients of free algebras by
homogeneous ideals, and the program GRB \cite{GRB} provides minimal
projective resolutions for finite dimensional modules over finite
dimensional quotients of path algebras.

The goal of this paper is to provide an algorithmic method that
constructs projective resolutions of modules over quotients of path
algebras.  We also modify this algorithm to construct minimal linear
projective resolutions of linear modules over a Koszul algebra in the
last section. The construction uses both the theory of projective
resolutions presented in \cite{GSZ} and the theory of Gr\"obner bases
for path algebras \cite{G2}. 

Before giving a brief summary of needed background material in the
next section, we provide an overview of the results of the paper. To
this end we recall the definition of a path algebra, with a fuller
explanation given in the next section.  Let $k$ be a field and let $Q$
be a finite quiver; that is, a finite directed graph. The path
algebra, $kQ$, has $k$-basis consisting of all the finite directed
paths in $Q$, and multiplication is induced by concatenation of paths.
For the remainder of this paper let $R=kQ$ denote a path algebra.  In
Section \ref{section:mainstep}, we present our main step in the
construction of a projective resolution of a module over a quotient of
a path algebra and then in Section \ref{section:construction}, we show
how to get a resolution using the main step. One of the more
interesting theoretical results is that if $I$ has a finite Gr\"obner
basis, and $M$ is a right $R/I$-module which is finitely presented as
an $R$-module, then the construction yields a projective
$R/I$-resolution of $M$ such that each projective occurring in the
resolution is finitely generated (even though $R/I$ need not be
noetherian). A discussion of the algorithmic aspects of the
construction of the resolution follows in the next section.  The final
section applies a modified version of our construction to give an
algorithmic method for constructing minimal linear projective
resolutions of linear modules over Koszul algebras. Applications of
this algorithm can be found in \cite{BGMS,BGSS}.

\section{Preliminaries}\label{section:preliminaries}
The main object of study in this paper is the construction of
projective resolutions of modules over quotients of path algebras. In
this section we recall notions and results on path algebras, Gr\"obner
basis theory for path algebras \cite{G2} and projective resolutions of
modules over quotients of these as presented in \cite{GSZ}.

Let $k$ be a field and let $Q$ be a finite quiver.  We denote the
vertex set of $Q$ by $Q_0$, the arrow set by $Q_1$ and let $\calB$
denote the set of finite directed paths in $Q$. The path algebra,
$kQ$, as a $k$-vector space, has basis $\calB$. Note that we view the
vertices of $Q$ as paths of length $0$.  If $p,q\in \calB$, then we
define $p\cdot q=(pq)$ if the terminus vertex of $p$ = the origin
vertex of $q$ and $0$ otherwise.  If $Q$ has a single vertex and $n$
arrows (loops), then $kQ$ is isomorphic to the free associative
algebra on $n$ noncommuting variables.  Hence the class of algebras we
study include quotients of free algebras. We refer the reader to
\cite{ARS} for a fuller description of path algebras and their
properties.

Beginning our background information, we summarize the theory of
projective resolutions presented in \cite{GSZ}. Let $I$ be an ideal in
$R=kQ$ and $\Lambda=R/I$.  Let $X$ be an $R$-$R$-bimodule (respectively a
left $R$-module, a right $R$-module) and $x$ in $X$. Then we say that
a nonzero element $x$ is \emph{uniform} (respectively \emph{left
uniform}, \emph{right uniform}) if there exist $u$ and $v$
(respectively $u$ in $Q_0$, $v$ in $Q_0$) in $Q_0$ such that $x=uxv$
(respectively $x=ux$, $x=xv$).  Note that if $Q$ has a single vertex
then every nonzero $x$ in $X$ is uniform (resp., left uniform, right
uniform).

If $X$ is a right $kQ$-module and $x$ is in $X$, then $\overline{x}$
denotes the natural residue class of $x$ in $X/XI$.

Suppose that $M$ is a finitely generated right $\Lambda$-module.
Then, as shown in \cite{GSZ}, there exist $t_n$ and $u_n$ in
$\{0,1,2,\dots\}\cup \infty$ with $u_0=0$, $\{f^n_i\}_{i\in
T_n=[1,\dots,t_n]}$, and $\{{f^n_i}'\}_{i\in U_n=[1,\dots,u_n]}$ such
that
\begin{enumerate}
\item[(i)] Each $f^0_i$ is a right uniform element of $R$ for all
$i\in T_0$.
\item[(ii)] Each $f_i^n$ is in $\amalg_{j\in
T_{n-1}}f^{n-1}_jR$ and is a right uniform element for all $i\in T_n$
  and all $n\ge 1$. 
\item[(iii)] Each ${f_i^n}'$ is in $\amalg_{j\in T_{n-1}}f^{n-1}_jI$ and
is a right uniform element for all $i\in U_n$ and all $n\ge 1$.
\item[(iv)] For each $n\ge 2$, \[(\amalg_{i\in T_{n-1}}f^{n-1}_iR)
\cap (\amalg_{i\in T_{n-2}}{f^{n-2}_i}I) =(\amalg_{i\in T_n} f^n_iR)
\amalg(\amalg_{i\in U_n}{f^n_i}'R).\]
\end{enumerate}

The next result explains how the above elements give rise to a
projective $\L$-resolution of $M$. For this we need some notation. Let
$f_1,\ldots,f_m$ be right uniform elements of $R$ and $v_1,\ldots,v_m$
vertices such that $f_iv_i=f_i$ for $i=1,\ldots,m$. For
$i=1,\ldots,m$, let $\varepsilon_i=(\varepsilon_{i1},\ldots,
\varepsilon_{im})$ in $\amalg_{i=1}^m f_iR$ be defined by
$\varepsilon_{ij}=0$ for $i\neq j$ and $\varepsilon_{ii}=f_i$. Let
$\overline{\varepsilon_i}$ in $\amalg_{i=1}^m f_iR/\amalg_{i=1}^m
f_iI$ be defined in a similar fashion.

\begin{thm}[\cite{GSZ}]\label{thm:fnresol} 
Let $M$ be a finitely generated right $\Lambda$-module and suppose
that, for $n \ge 0$, $t_n$ and $u_n$ are in $\{0,1,2,\dots\}\cup
\infty$, and $\{f^n_i\}_{i\in T_n=[1,\dots,t_n]}$, and \linebreak
$\{{f^n_i}'\}_{i\in U_n=[1,\dots,u_n]}$ are chosen satisfying
\emph{(i)-(iv)} above.  We have that $f^n_i=$ \linebreak $\sum_{j\in
T_{n-1}}f^{n-1}_jh^{n-1,n}_{j,i}$, for some right uniform elements
$h^{n-1,n}_{j,i}$ in $R$.  Let $L^n=(\amalg_{i\in
T_n}f^n_iR/\amalg_{i\in T_n}f^n_iI)$, and $e^{n+1}\colon
L^{n+1}\extto{} L^{n}$ be given by $\overline{f^n_jh^{n,n+1}_{j,i}}$
in the $j$-th component of $e^{n+1}(\overline{\varepsilon_i})$.  Then
\[\cdots \extto{e^{n+1}}L^n\extto{e^n}L^{n-1}\extto{e^{n-1}}\cdots \extto{e^1}
L^0\extto{} M\extto{} 0\]
is a projective $\Lambda$-resolution of $M$.
\end{thm}

A more precise statement of the goal of this paper is to show how to
algorithmically construct the $t_n$'s, the $u_n$'s, the $f^n_i$'s, and
the ${f^n_i}'$'s.  For this we need the theory of noncommutative
Gr\"obner bases in path algebras and we review this theory.  For more
complete details we refer the reader to \cite{G2}.

First we need to order the basis $\calB$ of paths in 
$Q$.  We say $>$ is an {\it admissible} order on $\calB$ if
the following properties hold.
\begin{enumerate}
\item The order $>$ is a well-order.
\item If $p,q\in\calB$ with $p>q$ then for all $r\in\calB$,
$pr>qr$ if both $pr$ and $qr$ are nonzero.
\item If $p,q\in\calB$ with $p>q$ then for all $r\in\calB$,
$rp>rq$ if both $rp$ and $rq$ are nonzero.
\item If $p=qr$ with $p,q$, and $r$ paths in $\calB$ then $p\ge q$ and
$p\ge r$.
\end{enumerate}

There are many admissible orders. For example, we arbitrarily order
the set of vertices of $Q$, and we arbitrarily order the set of arrows
of $Q$.  Set each vertex smaller that any arrow.  If $p$ and $q$ are
paths of length at least one, then $p>q$ if the length of $p$ is
greater than the length of $q$, or, the lengths are equal and
$p=a_1a_2\cdots a_n$ and $q=b_1b_2\cdots b_n$ with the $a_i$'s and
$b_i$'s arrows, then there is some $i$ such that $a_j=b_j$ if $j<i$
and $a_i>b_i$.  For the remainder of this section, let $>$ be an
admissible order on $\calB$.

If $x=\sum_{i=1}^n\alpha_ip_i\in kQ$ with $\alpha_i$ nonzero elements
of $k$ and $p_i$ distinct paths, then $\tip(x)=p_i$ if $p_i\ge p_j$
for $j=1,\dots, n$.  If $X\subseteq kQ$ then $\tip(X)=\{\tip(x) \mid
x\in X\setminus\{0\}\}$.  We say a subset $\calG$ of $I$ is a
\emph{Gr\"obner basis of $I$ (with respect to $>$)} if the ideal
generated by $\tip(\calG)$ equals the ideal generated by $\tip(I)$.

There is an extension of Buchberger's algorithm that allows one to
construct Gr\"obner bases for ideals.  A word of caution is needed
here, in that, in general, even if $I$ is finitely generated there may
not be a finite Gr\"obner basis for $I$.  But the ``algorithm''
sequentially constructs sets $\calG_1,\calG_2,\dots$ so that $\cup_i
\calG_i$ is a Gr\"obner basis and if there is a finite Gr\"obner basis
then the ``algorithm'' terminates in a finite number of steps.

We provide a small example of a Gr\"obner basis in our setting, which
we refer to later in the paper. 
\begin{example}\label{example:basic}
Let $Q$ be the quiver 
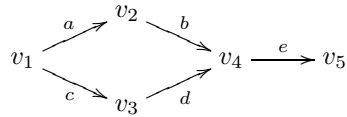
\[\xymatrix@R=6pt{ 
& v_2 \ar[dr]^b & & \\
v_1\ar[ur]^a\ar[dr]_c & & v_4\ar[r]^e & v_5\\
& v_3 \ar[ur]_d & & }\]
Let $I$ be the ideal generated by $ab-cd$ and $be$ in $kQ$. Choose $>_1$
to be the admissible order described earlier with
$v_5<v_4<v_3<v_2<v_1<e<d<c<b<a$. The algorithm described in \cite{G2}
yields the Gr\"obner basis $\calG=\{ab-cd,be,cde\}$.

If we change the order $>_1$ to $>_2$, where $v_5 < \ldots < v_1 < a <
b < c < d < e$, then one may check that the Gr\"obner basis now is
$\calG=\{ab-cd, be\}$ (since $\tip(ab-cd)=cd$).
\end{example}

We need the concept of a tip-reduced set of uniform elements of $R$.
If $p$ and $q$ are paths, we say {\it $p$ divides $q$}, denoted $p\mid
q$ (respectively, {\it $p$ left divides $q$} denoted $p\mid_lq$, and
{\it $p$ right divides $q$}, denoted $p\mid_rq$) if $q=rps$ for some
paths $r$ and $s$ (resp., if $q=ps$ for some path $s$, and $q=rp$ for
some path $r$).  We say a set of nonzero elements $X$ in $R$ is {\it
tip-reduced} if, for $x,y\in X$, $\tip(x)\mid \tip(y)$ implies $x=y$.
Since $>$ is a well-order on $\calB$, we have the following result.

\begin{prop}[\cite{G,G2}] \label{prop:tipred1}
If $X$ is a finite set of uniform elements of $R$, there is a finite
algorithm to produce a tip-reduced set of uniform elements $Y$ of $R$
such that the ideal generated by $X$ equals the ideal generated by
$Y$.
\end{prop}

To extend this concept to right projective $R$-modules, we must extend
the notion of a tip.  Let $\calI$ be an index set and, for each $i$,
let $v_i\in Q_0$.  Consider the right projective $R$-module
$P=\amalg_{i\in \calI} v_iR$.  Let $\calC$ be the set of all elements
of $P$ of the form $x=(x_i)_{i\in\calI}$ such that, for all but one
$i$, $x_i=0$, and, in that one coordinate, $x_i$ is a path (with
origin vertex $v_i$).  Then $\calC$ is a $k$-basis of $P$.  We now
define a well-order $>_P$ on $\calC$.  First let $>_{\calI}$ be a
well-order on $\calI$.  If $x=(x_i)$ and $y=(y_i)$ are elements of
$\calC$, then $x>_Py$ if the nonzero path occurring in $x$ is greater
than the nonzero path occurring in $y$ (using the admissible order $>$
on $\calB$), or, if these paths are equal, then the coordinate that
the nonzero entry occurs in $x$ is greater than the coordinate that
the nonzero entry occurs in $y$ (using $>_{\calI}$). The reader may
verify that $>_P$ is a well-order on $\calC$.  The order $>_P$ is
dependent on the choice of $>_{\calI}$ but in the remainder of this
paper, for each set $\calI$, we fix some well-order $>_{\calI}$.

Keeping the notation of the previous paragraph, if $w=(w_i)$ is a
nonzero right uniform element of $P$, we let $\tip(w)$ be the element
of $\calC$ such that the nonzero element of $\tip(w)$ is $p$ in
coordinate $i^*$ where (i) $\tip(w_{i^*})=p$, (ii) if $w_j\ne 0$ then
$p\ge \tip(w_j)$, and (iii) $i^*\geq _{\calI} j$ for all $j$ such that
$\tip(w_j)=p$.  We call $p$ the \emph{tip path of $w$} and denote it
by $\tippath(w)$, and we call $i^*$ the \emph{tip coordinate of $w$}
and denote it by $\tipcoord(w)$.  Letting
$\varepsilon_{i}=(\varepsilon_{ij})$ be in $P$ with
$\varepsilon_{ij}=\delta_{ij}v_i$, we see that
$\tip(w)=\varepsilon_{i^*}p$, where $p=\tippath(w)$ and
$i^*=\tipcoord(w)$.  The proof of the following result is left to the
reader after noting that if $x\in P$ and $p\in\calB$ such that
$\tip(x)p\ne 0$, then $\tip(xp)=\tip(x)p$,
$\tippath(xp)=\tippath(x)p$, and $\tipcoord(xp) =\tipcoord(x)$.

\begin{prop} Keeping the above notation, suppose $x$ and $y$ are right
uniform elements of $P$ and $p,q\in \calB$.  Then $\tip(x)>_P \tip(y)$
implies $\tip(xp)>_P\tip(yp)$, if both $\tip(x)p$ and $\tip(y)p$ are
nonzero.
\end{prop}

If $x,y\in \calC$, then we say $x$ {\it left divides} $y$, written as
$x\mid_ly$, if there is some path $p$ such that $xp=y$.  Note that
$x\mid_ly$ if and only if $\tippath(x)\mid_l\tippath(y)$ and
$\tipcoord(x)=\tipcoord(y)$.  We say a set of right uniform nonzero
elements $X$ of $P$ is {\it right tip-reduced} if, for each $x,y\in X$
with $\tip(x)\mid_l\tip(y)$ implies $x=y$.

\begin{prop}[\cite{G,G2}] \label{prop:tipred2} Let $P=\amalg_{i\in\calI}v_iR$
be as above.  If $X$ is a set of right uniform elements of $P$, then
there is a right tip-reduced subset $Y$ of $P$ of right uniform
elements such that the submodule of $P$ generated by $X$ equals the
submodule generated by $Y$.  Moreover, if $X$ is a finite set, then
there is a finite algorithm to produce such a $Y$ with $Y$ finite.
\end{prop}
\begin{proof}  
We include this proof for completeness. Let $A$ be the submodule
generated by $X$, and set $T=\{\tip(a)\mid a\in A\setminus
\{0\}\}$. Let $T^*=\{t\in T\mid \text{if\ } t'\mid_l t \text{\ and\ }
t'\in T \text{\ then\ } t=t'\}$. For each $t$ in $T^*$ choose a right
uniform $y_t$ in $A$ such that $\tip(y_t)=t$. Setting $Y=\{y_t\mid
t\in T^*\}$, we see that $Y$ is a right tip-reduced set. Let $B$
denote the submodule of $P$ generated $Y$. We claim that
$A=B$. Clearly $B\subseteq A$. Suppose that $A\not\subseteq B$, and
let $a$ in $A\setminus B$ such that $\tip(a)$ is minimal for
$\tip(a')$ for $a'$ in $A\setminus B$. Without loss of generality we
may suppose that $a$ is a right uniform element, since $a=\sum_{v\in
  Q_0} av$. Then, by definition of $T^*$, there is some $t$ in $T^*$
such that $\tip(t)\mid_l \tip(a)$. Then there is some path $p$ and $c$
in $k$ such that $a-cy_tp$ has smaller tip than $a$. But then
$a-cy_tp$ is in $B$. Since $cy_tp$ is in $B$, we have a
contradiction. 

Now suppose that $X$ is a finite set of right uniform elements in
$P$. Consider the following process. 
\begin{quote}

While $X$ is not tip-reduced, 
\begin{quote}
Let $X=\{x_1,\ldots,x_n\}$. Suppose $\tip(x_i)\mid_l \tip(x_j)$ for
some $i\neq j$. We let
$X_1=\{x_1,\ldots,x_{j-1},x'_j,x_{j+1},\ldots,x_n\}$, where
$x_j'=x_j-cx_ip$ for some $c$ in $k$ and $p$ in $\calB$ such that
$\tip(x_j')<\tip(x_j)$. We see that the right submodule of $P$
generated by $X_1$ is equal to $A$. Replace $X$ by $X_1$.
\end{quote}
Output: $X$.
\end{quote}
This process has to stop in a finite number of steps, since $>_P$ is a
well-order. This completes the proof.
\end{proof}

The importance of a generating set being right tip-reduced is
demonstrated by the following result.
\begin{prop}[\cite{G}]\label{prop:tip-reducedsum}
Let $P'$ be a submodule of a projective right $R$-module
$P=\amalg_{i\in \calI} v_iR$, where $\calI$ is an index set and $v_i$
is in $Q_0$ for all $i$. If $\{f_j\}_{j\in\mathcal{J}}$ is a right
tip-reduced generating set for $P'$ consisting of right uniform
elements, then $P'=\amalg_{j\in \mathcal{J}} f_jR$.
\end{prop}

We provide a small example to clarify some of ideas presented above.

\begin{example}
Let $Q$ be the quiver
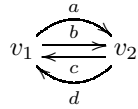
\[\xymatrix{v_1\ar@/^1pc/[r]^a\ar@<.5ex>[r]^b &
  v_2\ar@/^1pc/[l]^d\ar@<.5ex>[l]^c}\] 
with admissible order $>$ on $\calB$ defined earlier with
$v_1<v_2<a<b<c<d$. Let $P=v_1R\amalg v_1R\amalg v_2R\amalg v_2R$ with
order $>_P$ induced by $>$ and $(v_1,0,0,0)>(0,v_1,0,0)>(0,0,v_2,0)>
(0,0,0,v_2)$. Let $X=\{f_1,f_2,f_3\}$ in $P$, where
$f_1=(ac,bcac+acac,0,cac)$, $f_2=(v_1,bc+ac,d,0)$ and
$f_3=(a,b,da+cb,da)$. We right tip-reduce $X$. We see that
$\tip(f_1)=(0,bcac,0,0)$, $\tip(f_2)=(0,bc,0,0)$ and
$\tip(f_3)=(0,0,da,0)$. Since $\tip(f_2)\mid_l\tip(f_1)$, we replace
$f_1$ by $f'_1=f_1-f_2ac=(0,0,-dac,cac)$. Then
$\tip(f'_1)=(0,0,dac,0)$ and we see that
$\tip(f_3)\mid_l\tip(f'_1)$. Hence we replace $f'_1$ by
$f''_1=f'_1+f_3c=(ac,bc,abc,dac+cac)$. Then
$\tip(f''_1)=(0,0,0,dac)$. Hence the set $X^*=\{f''_1,f_2,f_3\}$ is
right tip-reduced and we have that (i) $X$ and $X^*$ generated the
same submodule, say $P'$ of $P$ and (ii) $P'=f''_1R\amalg f_2R\amalg
f_3R$.
\end{example}

We note that, given a set subset $X$ of $\amalg_{i\calI}v_iR$, there
is no unique right tip-reduced $X'$ generating the same submodule as
$X$. 

We need one final definition. Suppose $p$ and $q$ are paths. We say
\emph{$q$ and $p$ overlap} or \emph{$q$ overlaps $p$} if there exist
paths $r$ and $s$ such that $pr=sq$. We allow $s$ to be a vertex, but
$r$ is a path of length at least length $1$.  An overlap relation may
be illustrated in the following way,
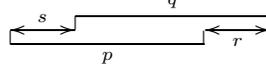
\[\xymatrix@W=0pt@M=0.3pt{%
\ar@{^{|}-^{|}}@<-1.25ex>[rrr]_p\ar@{{<}-{>}}[r]^{s} &
\ar@{_{|}-_{|}}@<1.25ex>[rrr]^{q} & & \ar@{{<}-{>}}[r]_r & }\] %

\section{The main step}\label{section:mainstep}

In this section we present the main step of our construction of a
projective resolution of a module over a quotient, $\L$, of a path
algebra. Beginning with a presentation of a $\L$-module over the
path algebra, we show how to find a presentation, over the path
algebra, of the first syzygy of the module over $\L$. This gives rise
to an inductive algorithm for finding a projective resolution of a
module over $\L$ described in the next section.

Let $I$ be an ideal in a path algebra $R=kQ$, let $\L=R/I$, and let
$\calG=\{ g^2_i\}_{i\in\calI}$ be a uniform, tip-reduced Gr\"obner
basis for the ideal $I$ with respect to some admissible order $>$.

Let $M$ be a right $\L$-module. By \cite{G} there exists an 
$R$-presentation of $M$ of the form
\[0\to (\amalg_{i\in T_1} f^1_iR)\amalg (\amalg_{j\in
U_1} {f^1_j}'R)\extto{H^1}\amalg_{i\in T_0} f^0_iR\extto{\pi} M\to 0,\]
where 
\begin{enumerate}
\item[(i)] $H^1$ is an inclusion,
\item[(ii)] $f^0_i$'s, $f^1_i$'s and ${f^1_i}'$ are right uniform, 
\item[(iii)] ${f^1_j}'$ is in $\amalg f^0_i I$ for all $j$ in $U_1$, 
\item[(iv)] the set $\{f^1_i\}_{i\in T_1}\cup \{{f^1_i}'\}_{i\in U_1}$ is
  right tip-reduced.
\end{enumerate}
Our goal is to construct sets $\{ f^2_i\}_{i\in T_2}$ and $\{
{f^2_i}'\}_{i\in U_2} $, such that $\{ f^2_i\}_{i\in T_2}\cup \{
{f^2_i}'\}_{i\in U_2} $ is a right uniform and right tip-reduced set in
$\amalg f^1_iR$, the set $\{ {f^2_i}'\}_{i\in U_2}$ is in $\amalg
f^1_i I$, and 
\[0\to (\amalg_{i\in T_2} f^2_iR)\amalg (\amalg_{j\in
U_2} {f^2_j}'R)\extto{H^2}\amalg_{i\in T_1} f^1_iR\to \Omega^1_\L(M)\to 0,\] 
is an exact sequence of right $R$-modules, where $H^2$ is an 
inclusion map and $\Omega^1_\L(M)$ is the kernel of $\amalg_{i\in T_0}
f^0_iR/\amalg_{i\in T_0} f^0_iI\to M$. 

Recall from \cite{GSZ} that we want to construct the $f^2_i$'s and the
${f^2_i}'$'s so that
\[(\amalg_{i\in T_1} f^1_iR)\cap (\amalg_{i\in
T_0} f^0_iI) = (\amalg_{i\in T_2} f^2_iR)\amalg (\amalg_{j\in
U_2} {f^2_j}'R).\]
This equality can be seen from the following short exact sequence of
right $R$-modules 
\[0\to (\amalg_{i\in T_1} f^1_iR)\cap (\amalg_{i\in T_0} f^0_iI)\to 
\amalg_{i\in T_1} f^1_iR \to \Omega^1_\L(M)\to 0\]
The existence of this exact sequence is obtained by considering the
exact sequence of right $R$-modules given by the left hand column of
the following commutative exact diagram 
\[\xymatrix{%
& 0\ar[d] & 0\ar[d] & & \\
& \amalg_{i\in T_0} f^0_iI \ar@{=}[r]\ar[d] & \amalg_{i\in T_0} f^0_iI
\ar[d] & & \\
0\ar[r] & (\amalg_{i\in T_1} f^1_iR)\amalg (\amalg_{j\in U_1}
{f^1_j}'R) \ar[r]\ar[d] & \amalg_{i\in T_0} f^0_iR \ar[r]^{\pi}\ar[d]
& M \ar[r]\ar@{=}[d] & 0\\
0\ar[r] & \Omega^1_\L(M)\ar[r]\ar[d] & 
\amalg_{i\in T_0} f^0_iR/\amalg_{i\in T_0} f^0_iI\ar[r]\ar[d] & M\ar[r] & 0\\
& 0 & 0 & & }\]

To construct the $f^2_i$'s, we need some preliminary definitions.  Let
$p$ be a path in $Q$ of length at least one.  We define $X(p)$ to be the
set of paths $q$ that satisfy the following conditions:
\begin{enumerate}
\item $p\mid_l q$. 
\item There is some $g^2_i\in \calG$ such that $\tip(g^2_i)\mid_r q$.  
\item If there are paths $r$ and $s$ and $g^2_j\in \calG$ such that
$q=r\tip(g^2_j)s$ then $s$ is a vertex (and hence
$i=j$ since $\{g^2_t\}_{t\in \calI}$ is tip-reduced).
\end{enumerate}
The following figures illustrate (1) and (2) in the definition of $X(p)$:
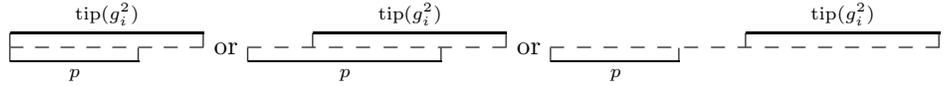
\[\xymatrix@W=0pt@M=0.3pt{%
\ar@{^{|}-^{|}}@<-1.25ex>[rr]_p\ar@{--}[rrr] 
\ar@{_{|}-_{|}}@<1.25ex>[rrr]^{\tip(g^2_i)} & & & }
\text{\ or\ }
\xymatrix@W=0pt@M=0.3pt{%
\ar@{^{|}-^{|}}@<-1.25ex>[rrr]_p\ar@{--}[rrrr] &
\ar@{_{|}-_{|}}@<1.25ex>[rrr]^{\tip(g^2_i)} & & & }
\text{\ or\ }
\xymatrix@W=0pt@M=0.3pt{%
\ar@{^{|}-^{|}}@<-1.25ex>[rr]_p\ar@{--}[rrrrrr] & & & 
\ar@{_{|}-_{|}}@<1.25ex>[rrr]^{\tip(g^2_i)}
& & & }
\] 
where $q$ is the path indicated by the dashed lines. 

If $q\in X(p)$ and $q=q'\tip(g^2_i)$, call $g^2_i$ the {\em end
relation of $q$}.  We break $X(p)$ into two disjoint sets.  Let
\[O(p)=\{q\in X(p) \mid \text{the tip of the end relation of }q\text{
and }p\text{ overlap}\}\] and
\[N(p)=X(p)\setminus O(p).\]
Elements $q$ in $O(p)$ can be describe by the following diagram
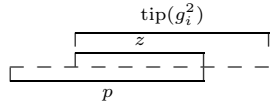
\[\xymatrix@W=0pt@M=0.3pt{%
\ar@{^{|}-^{|}}@<-1.25ex>[rrr]_p\ar@{--}[rrrr] & 
\ar@{_{|}-_{|}}@<3ex>[rrr]^{\tip(g^2_i)}\ar@{_{|}-_{|}}@<1.25ex>[rr]^z 
& & & }\] 
where $z$ is a path of length at least one (in particular we allow
$z=p$). Again, $q$ is the path indicated by the dashed line. Elements
$q$ in $N(p)$ are illustrated by the following diagram
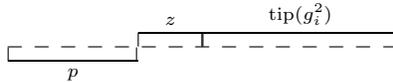
\[\xymatrix@W=0pt@M=0.3pt{%
\ar@{^{|}-^{|}}@<-1.25ex>[rr]_p\ar@{--}[rrrrrr] & &
\ar@{_{|}-_{|}}@<1.25ex>[r]^z  & 
\ar@{_{|}-_{|}}@<1.25ex>[rrr]^{\tip(g^2_i)}
& & & }\]
where $z$ is a path of length at least zero. 

We can now define $T_2$, the index set for the $f^2_i$'s, and $U_2$,
the index set for the ${f^2_i}'$'s.  Let $T_2 =\{(i,q) \mid {i\in
T_1}\text{ and }q\in O(\tippath(f^1_i))\}$ and $U_2 =\{(i,q) \mid
{i\in T_1}\text{ and }q\in N(\tippath(f^1_i))\}$.  We remark here that
$T_2$ and $U_2$ are countable sets, since $T_1$ and $\calB$ are
countable sets.  To define the $f^2_s$, suppose that $s=(i,q)\in T_2$
and that $\tipcoord(f^1_i)=i^*$. From the definition of $T_2$, we see
that $q=\tippath(f^1_i)p=q'\tip(g^2_j)$ for some paths $p$ and $q'$
and $g^2_j\in \calG$ is the end relation of $q$.  Consider
$f^1_ip-\varepsilon_{i^*}cq'g^2_j$, where $\varepsilon_{i^*}$ is
defined as in Section \ref{section:preliminaries} and
$c=\frac{\coefftip(f^1_i)}{\coefftip(g^2_j)}$ in $k$.  We see that
$f^1_ip-\varepsilon_{i^*}cq'g^2_j$ is right uniform. Note that
$\varepsilon_{i^*}cq'g^2_j$ is in $\amalg_{l\in T_0}f^0_lR$ and has
only one non-zero component, namely $cq'g^2_j$ in the same component
as $\tipcoord(f^1_i)$.  Clearly
$\pi(f^1_ip-\varepsilon_{i^*}cq'g^2_j)=0$, so that
$f^1_ip-\varepsilon_{i^*}cq'g^2_j$ is in $(\amalg_{i\in T_1}
f^1_iR)\amalg (\amalg_{j\in U_1} {f^1_j}'R)$. Hence,
\[f^1_ip-\varepsilon_{i^*}cq'g^2_j=\sum_{j\in T_1}f^1_jr_j+\sum_{j\in
  U_1} {f^1_j}'s_j\] for some $r_j$ and $s_j$ in $R$. By the unicity
of the sums, there is a vertex $v$ such that $f^1_ipv=f^1_ip$,
$\varepsilon_{i^*}cq'g^2_jv=\varepsilon_{i^*}cq'g^2_j$, $f^1_lr_lv=
f^1_lr_l$ for all $l$ in $T_1$ and ${f^1_l}'s_lv={f^1_l}'s_l$ for all
$l$ in $U_1$. Since $\{ f^1_l\}_{l\in T_1}\cup\{{f^1_l}'\}_{l\in U_1}$
is a right tip-reduced right Gr\"obner basis for $(\amalg_{l\in
T_1}f^1_lR)\amalg (\amalg_{l\in U_1} {f^1_l}'R)$ and since
$\tip(f^1_ip-\varepsilon_icq'g^2_j)< \tip(f^1_ip)$, we see that
$\tip(f^1_ip)>\tip(f^1_jr_j)$ for all $j\in T_1$. Let
$f^2_s=f^1_ip-\sum_{j\in T_1}f^1_jr_j$. Then, since
$\varepsilon_{i^*}cq'g^2_j$ and each ${f^1_j}'$ is in $\amalg_{u\in
T_0} f^0_uI$, the element $f^2_s$ is in $(\amalg_{j\in T_1}f^1_jR)\cap
(\amalg_{u\in T_0} f^0_uI)$. Moreover, we see that $f^2_s$ is right
uniform. Thus, for each $s$ in $T_2$, we have constructed an
$f^2_s$. Note that $\tip(f^2_s)=\tip(f^1_i)p$.

We now construct the ${f^2_l}'$'s.  Let $s=(i,q)\in U_2$.  From the
definition of $U_2$, there is a path $z$ and a $g^2_j\in\calG$ such
that $q=\tippath(f^1_i)z\tip(g^2_j)$.  Define
${f^2_s}'=f^1_izg^2_j$. We have that each ${f^2_s}'$ is in
$\amalg_{j\in T_1} f^1_jI$.  It is clear that ${f^2_s}'\in
(\amalg_{i\in T_1} f^1_iR)\amalg(\amalg_{j\in U_1} {f^1_j}'R)$ and
that $\tip({f^2_s}')=\tip(f^1_i)z\tip(g^2_j)$.

The next result proves the main properties of the $f^2_i$'s and the
${f^2_i}'$'s. 
\begin{thm}
\[(\amalg_{i\in T_1}f^1_iR)\cap (\amalg_{i\in T_0}f^0_iI) =
(\amalg_{i\in T_2} f^2_iR)\amalg(\amalg_{i\in U_2} {f^2_{i}}'R)\]
and 
\[\{f^2_i\}_{i\in T_2}\cup \{{f^2_{i}}'\}_{i\in U_2}\]
is right uniform and right tip-reduced and hence a right uniform and
right tip-reduced right Gr\"obner basis for $(\amalg_{i\in
T_1}f^1_iR)\cap (\amalg_{i\in T_0}f^0_iI)$. Furthermore, each
${f^2_s}'$ is in $\amalg_{j\in T_1} f^1_jI$. 
\end{thm}
\begin{proof}
We have seen that the $f^2_s$'s and the ${f^2_s}'$'s are in $(\amalg_{i\in
T_1}f^1_iR)\cap (\amalg_{i\in T_0}f^0_iI)$, that they are right
uniform elements, and that each ${f^2_s}'$ is in $\amalg_{j\in T_1}
f^1_jI$.

We note that if $s\in T_2$ with $s=(i,q)$ and $g^2_j$ is the end
relation of $q$, then there are paths $p$ and $q'$ such that
$\tippath(f^1_i)p=q'\tip(g^2_j)$. We have seen that
$\tippath(f^2_s)=\tippath(f^1_i)p$. Let $i^*=\tipcoord(f^1_i)$.  Then
$\tippath(f^1_i)p$ occurs in the $i^*$-th component of $f^2_s$ viewed
as an element of $\amalg_{l\in T_0}f^0_lR$.  On the other hand, if
$s=(i,q)\in U_2$ with $q$ having end relation $g^2_j$, then there is
path $z$ such that $q=\tippath(f^1_i)z\tip(g^2_j)$.  We see that
$\tippath{f^2_{s}}'=\tippath(f^1_i)z\tip(g^2_j)$ in the
$\tipcoord(f^1_{i})$ coordinate of $\amalg_{i\in T_0}f^0_lR$.

Next we show that $\{f^2_s\}_{s\in T_2}\cup \{{f^2_{s}}'\}_{s\in U_2}$
is right tip-reduced.  Suppose not.  Since $\{f^1_j\}_{j\in T_1}$ is
right tip-reduced, it is clear that $\{f^2_s\}_{s\in T_2}$ and
$\{{f^2_{s}}'\}_{s\in U_2}$ are both right tip-reduced sets.  Suppose
that, for some $s\in T_2$ and $s'\in U_2$,
$\tipcoord(f^2_s)=\tipcoord({f^2_{s'}}')$ and $\tippath(f^2_s)$ left
divides $\tippath({f^2_{s'}}')$.  Let $s=(i,q)$ and $s'=(i',q')$.  We
see that either $\tippath(f^1_i)$ left divides $\tippath(f^1_{i'})$ or
vise versa.  In either case, since $\{f^1_j\}$ is right tip-reduced,
we conclude that $i=i'$.  But then, since $\tippath(f^2_s)$ left
divides $\tippath({f^2_{s'}}')$, we conclude that the end relation of
$q$ appears before the end relation of $q'$ which contradicts property
(3) of the definition of $X(q')$.  Hence $\tippath(f^2_s)$ does not
left divides $\tippath({f^2_{s'}}')$.  A similar argument shows that
$\tippath({f^2_{s'}}')$ does not left divides $\tippath({f^2_s})$.  We
conclude that $\{f^2_s\}_{s\in T_2}\cup \{{f^2_{s}}'\}_{s\in U_2}$ is
right tip-reduced.

Since $\{f^2_s\}_{s\in T_2}\cup \{{f^2_{s}}'\}_{s\in U_2}$ is right
uniform right tip-reduced, the submodule generated by this set can be
written as $(\amalg_{s\in T_2}f^2_sR)\amalg(\amalg_{s\in
U_2}{f^2_{s}}'R)$ by Proposition \ref{prop:tip-reducedsum}.  It
remains to show $\{f^2_s\}_{s\in T_2}\cup \{{f^2_{s}}'\}_{s\in U_2}$
generates $(\amalg_{i\in T_1}f^1_iR)\cap (\amalg_{i\in T_0}f^0_iI)$.
We have already proven that
\[(\amalg_{s\in T_2}f^2_sR)\amalg(\amalg_{s\in U_2}{f^2_{s}}'R)\subseteq
(\amalg_{i\in T_1}f^1_iR)\cap (\amalg_{i\in T_0}f^0_iI).\] Suppose
that $\{f^2_s\}_{s\in T_2}\cup \{{f^2_{s}}'\}_{s\in U_2}$ does not
generate $(\amalg_{i\in T_1}f^1_iR)\cap (\amalg_{i\in
T_0}f^0_iI)$. Let $x\in (\amalg_{i\in T_1}f^1_iR)\cap (\amalg_{i\in
T_0}f^0_iI)$ such that $\tip(x)$ is minimal with respect to the
property that $x\notin (\amalg_{s\in T_2}f^2_sR)\amalg(\amalg_{s\in
U_2}{f^2_{s}}'R)$. Since $x$ is in $\amalg_{i\in T_1}f^1_iR$ and since
$f^1_i$'s are tip-reduced, it follows that $\tip(x)=\tip(f^1_i)p$ for
some $i$ in $T_1$ and some path $p$. On the other hand, $x$ is in
$\amalg_{i\in T_0}f^0_iI$, hence
$\tip(x)=\varepsilon_{i^*}q\tip(g^2_j)z$ for some $j$ in $T_0$ and
some paths $q$ and $z$. Thus
$\tippath(x)=\tippath(f^1_i)p=q\tip(g^2_j)z$. For all possible
$g^2_j$'s such that $\tip(x)=\varepsilon_{i^*}q\tip(g^2_j)z$, choose
$j$ such that $q$ has minimal length. Either $\tip(g^2_j)$ overlaps
$\tippath(f^1_i)$ or not.

If they do overlap, then there exists an $l$ in $T_2$ such that
$\tip(f^2_l)z=\tip(x)$. Since the tip of $x-cf^2_lz$ is smaller than
$\tip(x)$ for some $c$ in $k$, the difference $x-cf^2_lz$ is in
$(\amalg_{s\in T_2}f^2_sR)\amalg(\amalg_{s\in U_2}{f^2_{s}}'R)$. This
is a contradiction.

If $\tip(g^2_j)$ does not overlap $\tippath(f^1_i)$, then there is
some $l$ in $U_2$ such that $\tip({f^2_l}')z=\tip(x)$ for some path
$z$ in $Q$. A similar argument as above leads to a contradiction. This
completes the proof.
\end{proof}

\begin{example}\label{example:basicII}
We continue Example \ref{example:basic}. First we use the order
$>_1$. Let $\L=kQ/I$ where $I$ is generated by $ab-cd$ and $be$. Let
$M=v_1\L/\rrad$, where $\rrad$ is the Jacobson radical of $\L$. It is 
immediate that, $M$, as a right $R$-module, has a projective
presentation
\[0\to aR\amalg cR\extto{H^1} v_1R\to M\to 0,\]
where $H^1(a)=a$ and $H^1(c)=c$. Recall that
$\calG=\{ab-cd,be,cde\}$. Let $g^2_1=ab-cd$, $g^2_2=be$, and
$g^2_3=cde$. Then $\tip(g^2_1)=ab$, $\tip(g^2_2)=be$ and
$\tip(g^2_3)=cde$. Let $f^0_1=v_1$, $f^1_1=a$ and $f^1_2=c$. We see
that $T_1=\{1,2\}$. We find $T_2 =\{(i,q) \mid {i\in T_1}\text{ and
}q\in O(\tippath(f^1_i))\}$ and $U_2 =\{(i,q) \mid {i\in T_1}\text{
and }q\in N(\tippath(f^1_i))\}$. First note that
$X(\tippath(f^1_1))=X(a) =\{ab\}=O(a)$ and $X(\tippath(f^1_2))=X(c)
=\{cde\}=O(c)$. Hence $T_2=\{(1,ab),(2,cde)\}$ and
$U_2=\emptyset$. For $(1,ab)$ we calculate $f^1_1b-v_1g^2_1=
ab-v_1(ab-cd)=cd=f^1_2d$. Therefore
$f^2_{(1,ab)}=f^1_1b-f^1_2d=ab$. Similarly we see that
$f^2_{(2,cde)}=f^1_2de=cde$. Thus we have
\[f^2_{(1,ab)}R\amalg f^2_{(2,cde)}R\extto{H^2} f^1_1R\amalg f^1_2R\]
with $H^2(f^2_{(1,ab)})=f^1_1b-f^1_2d$ and
$H^2(f^2_{(2,cde)})=f^1_2de$.

If we change the order to $>_2$, we still have $f^0_1=v_1$, $f^1_1=a$
and $f^1_2=c$. But now $\calG=\{ab-cd,be\}$ with $\tip(ab-cd)=cd$ and
$\tip(be)=be$. The reader may check that $T_2=\{(2,cd)\}$,
$U_2=\emptyset$, $f^2_{(2,cd)}=ab-cd$ and 
\[f^2_{(2,cd)}R\extto{H^2} f^1_1R\amalg f^1_2R,\]
where $H^2(f^2_{(2,cd)})=f^1_1b-f^1_2d$.
\end{example}

\section{Constructing resolutions}\label{section:construction}

This section is devoted to constructing a projective resolution of a
module over a quotient of a path algebra using the main step of the
previous section.

Let $M$ be a right $\L$-module and suppose that we have an
$R$-presentation of $M$ of the form
\begin{equation}\label{eq:Rpresentation}
0\to (\amalg_{i\in T_1} f^1_iR)\amalg (\amalg_{j\in
U_1} {f^1_j}'R)\extto{H^1}\amalg_{i\in T_0} f^0_iR\extto{\pi} M\to
0,
\end{equation}
where the $f^0_i$'s, $f^1_i$'s and ${f^1_i}'$'s are right uniform
elements, each ${f^1_i}'$ is in $\amalg_{i\in T_0} f^0_iI$, and 
$\{f^1_i\}_{i\in T_1}\cup \{{f^1_i}'\}_{i\in U_1}$ is right
tip-reduced. 

In the previous section we showed how to construct $f^2_i$'s and
${f^2_i}'$'s such that 
\[(\amalg_{i\in T_1}f^1_iR)\cap (\amalg_{i\in T_0}f^0_iI) =
(\amalg_{i\in T_2} f^2_iR)\amalg(\amalg_{i\in U_2} {f^2_{i}}'R)\]
and 
\[0\to (\amalg_{i\in T_2} f^2_iR)\amalg(\amalg_{i\in U_2}
       {f^2_{i}}'R)\extto{H^2}  
\amalg_{i\in T_1} f^1_iR \to \Omega^1_\L(M)\to 0\]
is an exact sequence of right $R$-modules, where $\Omega^1_\L(M)$ is
$\Ker(\amalg_{i=1}^{t_0}f^0_iR/\amalg_{i=1}^{t_0}f^0_iI\to M)$. 

From our construction, (i) $H^2$ is an inclusion map, (ii) the
elements $f^2_i$'s and ${f^2_i}'$'s are right uniform, (iii) each
${f^2_i}'$ is in $\amalg_{i\in T_1}f^1_iI$, and (iv) the set
$\{f^2_i\}_{i\in T_2}\cup \{{f^2_i}'\}_{i\in U_2}$ is right
tip-reduced.  Replacing $M$ by $\Omega_\L^1(M)$, we may view $f^1_i$'s
as $f^0$'s and the $f^2_i$'s as $f^1$'s and, applying our main step,
we may construct elements $f^3_i$'s and ${f^3_i}'$'s in
$\amalg_{i=1}^{t_2}f^2_iR$ so that
\[(\amalg_{i\in T_2}f^2_iR)\cap (\amalg_{i\in T_1}f^1_iI) =
(\amalg_{i\in T_3} f^3_iR)\amalg(\amalg_{i\in U_3} {f^3_{i}}'R)\]
and 
\[0\to (\amalg_{i\in T_3} f^3_iR)\amalg(\amalg_{i\in U_3}
       {f^3_{i}}'R)\extto{H^3} \amalg_{i\in T_2} f^2_iR \to
\Omega^2_\L(M)\to 0\] 
is an exact sequence of right $R$-modules with $\Omega^2_\L(M)$ being 
$\Ker(\amalg_{i=1}^{t_1}f^1_iR/\amalg_{i=1}^{t_1}f^1_iI\to
\Omega_\L^1(M))$, (i) $H^3$ is an inclusion
map, (ii) the elements $f^3_i$'s and ${f^3_i}'$'s are right uniform,
(iii) each ${f^3_i}'$ is in $\amalg_{i\in T_2}f^2_iI$, and (iv) the
set $\{f^3_i\}_{i\in T_3}\cup \{{f^3_i}'\}_{i\in U_3}$ is right
tip-reduced.

Repeating the above procedure, we obtain, for $n\geq 2$ elements
$f^n_i$'s and ${f^n_i}'$'s in $\amalg_{i=1}^{t_{n-1}} f^{n-1}_iR$ so
that
\[(\amalg_{i\in T_{n-1}}f^{n-1}_iR)\cap (\amalg_{i\in T_{n-2}}f^{n-2}_iI) =
(\amalg_{i\in T_n} f^n_iR)\amalg(\amalg_{i\in U_n} {f^n_{i}}'R)\]
and 
\[0\to (\amalg_{i\in T_n} f^n_iR)\amalg(\amalg_{i\in U_n}
       {f^n_{i}}'R)\extto{H^n}  
\amalg_{i\in T_{n-1}} f^{n-1}_iR \to \Omega^{n-1}_\L(M)\to 0\]
is an exact sequence of right $R$-modules with $\Omega^{n-1}_\L(M)$ being 
$\Ker(\amalg_{i=1}^{t_{n-2}}f^{n-2}_iR/\amalg_{i=1}^{t_{n-2}}f^{n-2}_iI\to
\Omega_\L^{n-2}(M))$, (i) $H^n$ is an inclusion
map, (ii) the elements $f^n_i$'s and ${f^n_i}'$'s are right uniform,
(iii) each ${f^n_i}'$ is in $\amalg_{i\in T_{n-1}}f^{n-1}_iI$, and
(iv) the set $\{f^n_i\}_{i\in T_n}\cup \{{f^n_i}'\}_{i\in U_n}$ is
right tip-reduced.

Since each $f^n_i$ is in $\amalg_{l\in T_{n-1}} f^{n-1}_lR$, we may
write 
\[f^n_i = \sum_{l\in T_{n-1}} f^{n-1}_lh^{n-1,n}_{li}\]
for some elements $h^{n-1,n}_{li}$ in $R$. Let 
\[L^n= \amalg_{l\in T_n} f^n_lR/\amalg_{l\in T_n} f^n_lI.\]
Let $v^n_i$ be the vertex in $Q$ such that $f^n_iv^n_i=f^n_i$. We see
that $L^n$ is isomorphic to $\amalg_{i\in T_n} v^n_i\Lambda$ for all
$n\geq 0$, hence it is a projective $\L$-module. Define $e^{n+1}\colon
L^{n+1}\to L^n$ by $e^{n+1}(\overline{f^{n+1}_i})$ equals
$\overline{f^n_jh^{n,n+1}_{ji}}$ in the component of $L^n$ corresponding to
$\overline{f^{n}_j}$. Now applying Theorem \ref{thm:fnresol} we
conclude that the resolution $(\mathbb{L},e)$
\[\cdots \extto{e^{n+1}}L^n\extto{e^n}L^{n-1}\extto{e^{n-1}}\cdots \extto{e^1}
L^0\extto{} M\extto{} 0\] 
is a projective $\L$-resolution of $M$. We call $(\mathbb{L},e)$
\emph{the resolution associated to \eqref{eq:Rpresentation}}.

\begin{example}\label{example:basicIII}
We now continue Example \ref{example:basicII}. Under the ordering
$>_1$ we have $f^0_1=v_1$, $f^1_1=a$, $f^1_2=c$,
$f^2_1=f^1_1b-f^1_2d$, and $f^2_2=f^1_2de$. We find the
$f^3_i$'s. Write $T_2$ as $\{1,2\}$. Then $\tip(f^2_1)=(b,0)$ and
$\tip(f^2_2)=(0,de)$ in $f^1_1R\amalg f^1_2R$. Hence
$X(\tippath(f^2_1))=X(b)=\{be\}$ and
$X(\tippath(f^2_2))=X(de)=\emptyset$. Thus $T_3=\{(1,be)\}$ and
$U_3=\emptyset$. For $(1,be)$, we calculate
\[f^2_1e-f^1_1g^2_2=f^1_1be-f^1_2de-f^1_1be=-f^1_2de=-f^2_2.\]
Hence $f^3_1=f^2_1e+f^2_2$ and we get 
\[f^3_1R\extto{H^3} f^2_1R\amalg f^2_2R,\]
where $H^3(f^3_1)=(e,v_5)$. The reader may check that
$T_4=\emptyset=U_4$. The induced resolution for $M$ over $\L$ by our
algorithm is 
\[0\to v_5\L\extto{\left(\begin{smallmatrix} e \\ v_5\end{smallmatrix}\right)}
v_4\L\amalg v_5\L \extto{\left(\begin{smallmatrix} b & 0\\ -d &
    de\end{smallmatrix}\right)} v_2\L\amalg
    v_3\L\extto{\left(\begin{smallmatrix} a 
    & c\end{smallmatrix}\right)} v_1\L\to M\to 0,\] 
since, for example, $f^3_1R/f^3_1I\simeq v_5\L$ and
    $f^2_1R/f^2_1I\simeq v_2\L$. 

For the order $>_2$ the reader may check that $T_3=\emptyset=U_3$ and
the induced resolution for $M$ over $\L$ is 
\[0\to v_4\L \extto{\left(\begin{smallmatrix} b\\
    -d\end{smallmatrix}\right)} v_2\L\amalg 
    v_3\L\extto{\left(\begin{smallmatrix} a 
    & c\end{smallmatrix}\right)} v_1\L\to M\to 0.\] 
\end{example}
We note that the resolution for the ordering $>_2$ is minimal whereas
the resolution for the ordering $>_1$ is not minimal. This example
shows that the constructed resolution is dependent on the choice of the
admissible order, since both the Gr\"obner basis for $I$ and the tips
are order dependent. An algorithmic method for minimizing a
non-minimal projective resolution of a finite dimensional module over
a finite dimensional quotient of a path algebra, is given in
\cite{G3,GSZ}. 

In the next section we discuss some algorithmic aspects of the above
construction. We mention that a special case of the results can be
found in \cite{A,AG}, where it is shown that simple modules of the
form $vR/J$, where $J$ is the ideal in $R$ generated by the arrows of
$Q$ and $v$ is a vertex.

We end this section by providing sufficient conditions for the
constructed resolution to have the property that each $L^n$ is
finitely generated.

\begin{prop}\label{prop:f1finite}
Let $\L=R/I$, where $R=kQ$ for some quiver $Q$. Suppose that there is
an admissible order $>$ on $\calB$ such that the Gr\"obner basis for
$I$ with respect to $>$ is finite.  Let $M$ be a right $\L$-module,
which, as a right $R$-module, has a presentation
\begin{equation}\label{eq:f1finite}
0\to (\amalg_{i\in T_1} f^1_iR)\amalg (\amalg_{j\in
U_1} {f^1_j}'R)\extto{H^1}\amalg_{i\in T_0} f^0_iR\extto{\pi} M\to
0
\end{equation}
with $T_0$ and $T_1$ finite sets, where $H^1$ is an inclusion,
$\{f^1_i\}_{i\in T_1}$ and $\{{f^1_i}'\}_{i\in U_1}$ are right uniform
and right tip-reduced sets, and ${f^1_i}'$'s are in $\amalg_{i\in
  T_0}f^0_iI$. Then there is a projective resolution $(\mathbb{L},e)$
of $M$ as a $\L$-module associated to \eqref{eq:f1finite} with the
property that each $L^n$ is finitely generated. 
\end{prop}
\begin{proof}
Let $\calG$ be a finite Gr\"obner basis of uniform elements for $I$.
Since $T_1$ and $\calG$ are finite sets, it follows that for each
$f^1_i$ there are only a finite number of $g^2_j$ such that
$\tippath(f^1_i)$ overlaps $\tip(g^2_j)$. It follows that $T_2$ is
also a finite set. Inductively we conclude that each $T_n$ is a finite
set for all $n\geq 0$.
\end{proof}
Note that in the previous result the set $U_1$ can be infinite. The
next result shows that if $M$ is finitely presented as a right
$R$-module, then all the sets $T_0$, $T_1$ and $U_1$ can be chosen to
finite in \eqref{eq:f1finite}. 
\begin{prop}\label{prop:fingenresol}
Let $\L=R/I$, where $R=kQ$ for some quiver $Q$.  Let $M$ be a right
$\L$-module which, as a right $R$-module, is finitely
presented. Suppose that there is an admissible order $>$ on $\calB$
such that the Gr\"obner basis for $I$ with respect to $>$ is
finite. Then there is a presentation of the form
\eqref{eq:Rpresentation} such that the the resolution $(\mathbb{L},e)$
associated to \eqref{eq:Rpresentation} has the property that each
$L^n$ is a finitely generated $\L$-module. The claim follows from
this. 
\end{prop}
\begin{proof}
Every projective right $R$-module is of the form $\amalg_{i\in\calI}
v_iR$, where $\calI$ is an index set and each $v_i$ is a vertex in
$Q_0$ \cite{G}.  Since $R$ is a hereditary algebra and since $M$ is a
finitely presented right $R$-module, we have a presentation of the
form
\[0\to \amalg_{i=1}^{n_1} w_iR\extto{\varphi} \amalg_{i=1}^{n_0}
v_iR\to M\to 0,\] 
where each $v_i$ and $w_i$ are vertices in $Q_0$. Let
$h^i=\varphi(w_i)$, which is a right uniform element for all
$i=1,\ldots, n_1$. Right tip-reduce the set $\{h^1,\ldots,h^{n_1}\}$,
and break the elements into two sets $\{f^1_1,\ldots,f^1_{t_1}\}$ and
$\{{f^1_1}',\ldots, {f^1_{u_1}}'\}$ so that each ${f^1_j}'$ is in
$\amalg_{i=1}^{n_0} v_iI$.  Finally set $t_0=n_0$ and $f^0_i=v_i$ for
$i=1,\ldots,t_0$. Thus we obtain the following presentation
\[0\to (\amalg_{i\in T_1} f^1_iR)\amalg (\amalg_{j\in
U_1} {f^1_j}'R)\extto{H^1}\amalg_{i\in T_0} f^0_iR\extto{\pi} M\to
0\]
of $M$ as a right $R$-module, where $f^0_i$'s, $f^1_i$'s, and
${f^1_i}'$'s are right uniform elements, both $T_0$ and $T_1$ are
finite sets, and the set $\{f^1_i\}_{i\in T_1}\cup \{{f^1_i}'\}_{i\in
U_1}$ is right tip-reduced.

We now apply Proposition \ref{prop:f1finite} to obtain our desired
result. 
\end{proof}

The previous result raises the question: Which right $\L$-modules are
finitely presented as right $R$-modules? The next result shows that
all finite dimensional right $\L$-modules are finitely presented as
right $R$-modules. 

\begin{prop}
Let $\L=R/I$, where $R=kQ$ for some quiver $Q$.  Let $M$ be a finite
dimensional right $\L$-module. Then $M$, as a right $R$-module, is
finitely presented. Furthermore, if $I$ has a finite Gr\"obner basis,
then there is a presentation of the form \eqref{eq:Rpresentation} such
that the resolution $(\mathbb{L},e)$ associated to
\eqref{eq:Rpresentation} has the property that each $L^n$ is a
finitely generated $\L$-module.
\end{prop}
\begin{proof} Let $M$ be a finite dimensional right $\L$-module.  It
is enough to show that $M$ is a finitely presented right
$R$-module. Let $A$ be the right annihilator of $M$ as a right
$R$-module, and let $\Gamma=kQ/A$. Then the $k$-algebra $\Gamma$ is
finite dimensional, and $M$ is a finitely generated right
$\Gamma$-module. Let $\{m_i\}_{i=1}^{t_0}$ be a finite set of right
uniform generators for $M$ as a $\Gamma$-module, and suppose that
$\{f^0_i\}_{i=1}^{t_0}$ is a set of vertices in $Q$ such that
$m_if^0_i=m_i$ for all $i=1,\ldots,t_0$. Since $\Gamma$ and $M$ are
finite dimensional, there is a projective $\Gamma$-presentation 
\[\amalg_{i=1}^dw_i\Gamma \to \amalg_{i=1}^{t_0} f^0_i\Gamma\to M\to
0,\] 
for some vertices $w_i$ in $Q$ and for some positive integer $d$. We
also have an exact sequence of right $R$-modules
\[0\to K\to \amalg_{i=1}^{t_0} f^0_iR\to M\to 0.\]
It can be seen that $(\amalg_{i=1}^dw_iR)\amalg (\amalg_{i=1}^{t_0}
f^0_iA)$ maps onto $K$. To show that $K$ is finitely generated as an
$R$-module, we need to show that $\amalg_{i=1}^{t_0} f^0_iA$ is
finitely generated. By \cite{G2} $A$ has a finite Gr\"obner basis with
respect to any admissible order. From \cite[Proposition 7.1]{G} and
the fact that $\Gamma$ is finite dimensional, it follows that
$\amalg_{i=1}^{t_0} f^0_iA$ is finitely generated. This shows that $M$
is a finitely presented right $R$-module. The final statement follows
from Proposition \ref{prop:fingenresol}. 
\end{proof}

For finite dimensional algebras $\L=kQ/I$, we have the following
consequence, since $I$ has a finite Gr\"obner basis with respect to
any admissible order \cite{G2}.
\begin{cor}
Let $\L=kQ/I$ be a finite dimensional algebra. Then any finitely
generated right $\L$-module has a projective $\L$-resolution
$(\mathbb{L},e)$ which can be constructed algorithmically such that
each $L^n$ is finitely generated.\qed
\end{cor}

\section{Algorithmic aspects}

In this section we discuss computational questions related to the
construction presented in the previous sections.  Our goal is to
clarify when we have actual (finite) algorithms for constructing
projective resolutions of modules over quotients of path algebras and
to provide an overview of the algorithms needed.  More precisely, let
$Q$ be a quiver, $I$ an ideal in $R=kQ$, and $\Lambda=R/I$.  Suppose
$M$ is a right $\Lambda$-module.  We wish to find conditions so that,
given a positive integer $N$, there is an algorithm based on the
construction in the earlier sections whose output is a projective
$\Lambda$-resolution
\[L^N\extto{e^N}L^{N-1}\extto{}\cdots \extto{}L^0\extto{}M\extto{}0.\]
We also discuss the input for such an algorithm.

There are two conditions needed; one on the ideal $I$ and one on the
module $M$.  We begin with the condition on the ideal $I$.  Let $>$ be
an admissible order on $\calB$, and $\calG$ a tip-reduced Gr\"obner
basis for $I$ with respect to $>$ consisting of uniform elements.  The
construction of $\calG$, given a finite set of generators of $I$ is
discussed in \cite{G2}.  For there to be a finite algorithm for
constructing $\calG$, we must assume that $\calG$ is finite.  As noted
earlier, if $R/I$ is finite dimensional over $k$, then there is finite
tip-reduced Gr\"obner basis for $I$ consisting of uniform elements. We
actually need something stronger than the existence of a finite
Gr\"obner basis.

Let 
\begin{multline}
\rtG = \{pg \mid p\notin \tip(I),  g\in \calG, 
\mbox{ if }r\tip(g')s=p\tip(g), \notag \\
\mbox{ for some }g'\in \calG, r,s\in \calB,
\mbox{ then } s\in Q_0 \mbox{ and }g=g'\}
\end{multline}
It is shown in \cite{G} that $\rtG$ is a right uniform, right
tip-reduced, right Gr\"obner basis for $I$.  If $\rtG$ is infinite, we
will not in general have a finitely terminating algorithm to right
tip-reduce sets needed in the construction.  For this reason, we need
to assume that $\rtG$ is finite. Of course, $\rtG$ being finite
implies that $\calG$ is a finite set.  We hasten to add that $\rtG$ is
finite if $R/I$ is finite dimensional over $k$, with $|\rtG|\le
\dim_k(\Lambda)\cdot |\calG|$.

We now consider the class of modules for which we have an algorithm to
construct a projective resolution.  Let $\Lambda=R/I$ and let $M$ be a
right $\Lambda$-module.  Since $\sum_{v\in Q_0}v=1$, we see that $M$
has a projective presentation as a right $\Lambda$-module of the form
\begin{equation}\label{eq:Lpresentation}
\amalg_{i\in\calI}w_i\Lambda \extto{\varphi}
\amalg_{i\in\calI'}v_i\Lambda\to M\to 0,
\end{equation}
where $\calI$ and $\calI'$ are index sets and each $v_i$ and $w_i$ are
vertices.  The assumption on $M$ that we need is that the index sets
$\calI$ and $\calI'$ in (\ref{eq:Lpresentation}) are finite.  The next
result is fundamental to the existence of an algorithm.

\begin{prop}\label{prop:liftpresentation}  Let $\Lambda=R/I$ where
$R=kQ$ for some quiver $Q$.  Assume that $\calG$ is a right
tip-reduced, right Gr\"obner basis for $I$ consisting of uniform
elements and assume further that $\rtG$ is finite.  Let $M$ be a right
$\Lambda$-module such that $M$ has a projective presentation as a
$\Lambda$-module of the form \eqref{eq:Lpresentation} with $\calI$ and
$\calI'$ finite sets.  Then, there is an algorithm, whose input is
\eqref{eq:Lpresentation} and output is nonnegative integers $t_0$,
$t_1$, and $u_1$ and a projective presentation of $M$ as an $R$-module
\[0\to(\amalg_{i=1}^{t_1}f^1_iR)\amalg(\amalg_{i=1}^{u_1}{f^1_i}'R)\extto{H^1}
\amalg_{i=1}^{t_0}f^0_iR\to M\to 0,\]
where \begin{enumerate}
\item $H^1$ is an inclusion map,
\item the $f^0_i$'s, $f^1_i$'s and ${f^1_i}'$'s are right uniform elements,
\item ${f^1_i}'\in \amalg_{i=1}^{t_0}f^0_iI$, for all $i=1,\dots,u_1$, and
\item $\{f^1_i\}_{i=1}^{t_1}\cup \{{f^1_i}'\}_{i=1}^{u_1}$ is right
tip-reduced.
\end{enumerate}
\end{prop}
\begin{proof}  By hypothesis, there exist nonnegative integers $t_0$ and $d$, 
vertices $w_i$, for $i=1,\dots,d$, and vertices $v_j$, for
$j=1,\dots,t_0$ such that there is an exact sequence of right
$\Lambda$-modules
\[\amalg_{i=1}^dw_i\Lambda \extto{\varphi}
\amalg_{j=1}^{t_0}v_i\Lambda\extto{\pi} M\to 0.\] For $j=1,\dots,t_0$,
let $f^0_j=v_j$.  We see that the surjection
$\pi\colon\amalg_{j=1}^{t_0}f^0_i\Lambda\to M$ induces a surjection
$\psi\colon\amalg_{j=1}^{t_0}f^0_iR\to M$ with kernel $K$, where
$\psi(f^0_i)=\pi(f^0_i)$.  It follows that there is a surjection
$\mu\colon (\amalg_{i=1}^{d}w_iR) \amalg(\amalg_{j=1}^{t_0}f^0_iI) \to
K$ since the kernel of $\amalg_{j=1}^{t_0}f^0_iR\to
\amalg_{j=1}^{t_0}f^0_i\L$ is $\amalg_{j=1}^{t_0}f^0_iI$.  The
surjection $\mu$ can be obtained algorithmically as follows.  For each
$i=1,\dots,d$, let
$\varphi(w_i)=x_i=(x_{i,1},\dots,x_{i,t_0})\in\amalg_{j=1}^{t_0}f^0_i\Lambda$.
Note that, $f^0_ix_{i,l}=x_{i,l}$ is in $f^0_i\Lambda$ for each $i$
and $l$.  Choose right uniform elements $h_{i,l}\in R$ such that
$\overline{h_{i,l}}=x_{i,l}$ and $h_{i,l}=f^0_ih_{i,l}$.  (Of course,
computationally, using Gr\"obner basis theory, one is representing the
$x_{i,l}$ as some $h_{i,l}$ already!) Let
$h_i=(h_{i,1},\dots,h_{i,t_0})$ in $\amalg_{i=1}^{t_0}f^0_iR$.  By our
assumptions, each $f^0_iI$ is a finitely generated right $R$-module
since the nonzero elements in $f^0_i\rtG$ form a right Gr\"obner basis
of $f^0_iI$.

We have that $\{h_{i}\}_{i=1}^d\cup \{f^0_i\rtG\}_{i=1}^{t_0}$ is a
finite right uniform generating set for $K$.  Let
$\{f^*_i\}_{i=1}^{d^*}$ be a right tip-reduced, right uniform set
obtained by right tip-reducing the set $\{h_i\}_{i=1}^d\cup
\{f^0_i\rtG\}_{i=1}^{t_0}$. Then $\{f^*_i\}_{i=1}^{d^*}$ is a right
tip-reduced and right uniform generating set for $K$ and
$K=\amalg_{i=1}^{d^*}f^*_iR$ by Proposition \ref{prop:tip-reducedsum}.
Since right tip-reduction is algorithmic, and right tip-reduction of a
right uniform set remains right uniform, taking the $f^1_i$'s to be
those $f^*_i$'s not in $\amalg_{j=1}^{t_0}f^0_jI$ and the ${f^1_i}'$'s
to be those $f^*_i$'s in $\amalg_{j=1}^{t_0}f^0_jI$, the result
follows.
\end{proof}

For the remainder of this section, we let $\Lambda=R/I$ where $R=kQ$
for some quiver $Q$ and assume that $\calG$ is a tip-reduced Gr\"obner
basis for $I$ consisting of uniform elements. Let $M$ be a right
$\Lambda$-module.  We keep the following two assumptions.  First, we
assume that $\rtG$ is finite.  Second, we assume that $M$ has a
projective presentation as a $\Lambda$-module of the form
(\ref{eq:Lpresentation}) with $\calI$ and $\calI'$ finite sets.

By Proposition \ref{prop:liftpresentation}, there is an algorithm, which we 
call {\bf LiftPresentation}, whose
input is a projective $\Lambda$-presentation of $M$ of form
(\ref{eq:Lpresentation})  
with $\calI$ and $\calI'$ finite, and whose output is nonnegative integers
$t_0$, $t_1$, and $u_1$ and a projective presentation of $M$ as
an $R$-module 
\[0\to(\amalg_{i=1}^{t_1}f^1_iR)\amalg(\amalg_{i=1}^{u_1}{f^1_i}'R)\extto{H^1}
\amalg_{i=1}^{t_0}f^0_iR\to M\to 0,\]
where \begin{enumerate}
\item $H^1$ is an inclusion map,
\item the $f^0_i$'s, $f^1_i$'s and ${f^1_i}'$'s are right uniform elements,
\item ${f^1_i}'\in \amalg_{i=1}^{t_0}f^0_iI$, for all 
  $i=1,\dots,u_1$, and
\item $\{f^1_iR)\}_{i=1}^{t_1}\cup \{{f^1_i}'\}_{i=1}^{u_1}$ is right
tip-reduced.
\end{enumerate}

Let $T$ be some finite set, for $i\in T$, let $\{f_i\}_{i\in T}$ be a
set of right uniform elements in $R$.  If $h_1,\dots, h_m,
{h_1}',\dots,{h_n}'$ is a right tip-reduced, right uniform subset of
$\amalg_{i\in T}f_iR$ and $x\in (\amalg_{i=1}^mh_iR)\amalg
(\amalg_{i=1}^n {h_i}'R)$ is right uniform, let {\bf FirstPart} be the
algorithm that takes as input $x$, $\{h_1,\dots, h_m\}$, and
$\{{h_1}',\dots,{h_n}'\}$ and outputs $\sum_{i=1}^mh_ir_i$ where
$x=(h_1r_1,\dots,h_mr_m,{h_1}'s_1,\dots,{h_n}'s_n)$ where the $r_i$
and $s_i$ are uniform elements of $R$. Note that {\bf FirstPart} is an
algorithm, since the $r_i$'s and the $s_i$'s can be obtained by right
tip-reducing $x$ by $\{h_1,\dots, h_m, {h_1}',\ldots,{h_n}'\}$.

If $\{h_1,\dots,h_m\}$ is a right uniform, right tip-reduced subset of
$\amalg_{i=1}^nf_iR$ where $\{f_i\}$ is a right uniform, right
tip-reduced set, let {\bf CreateMatrix} be the algorithm with input
$\{h_1,\dots,h_m\}$ and $\{f_1,\dots,f_n\}$ and output the $n\times m$
matrix $(h_{i,j})$ with uniform entries given by
$h_j=(f_1h_{1,j},\dots,f_nh_{n,j})$. Note that in {\bf CreateMatrix}
writing $h_j= (f_1h_{1,j},\dots,f_nh_{n,j})$ can be done
algorithmically by right tip-reducing $h_j$ by the set
$\{f_i\}_{i=1}^n$.

We now give an algorithmic description of the construction of a
projective resolution of a module $M$ given is the preceding sections.
We are given a field $k$, quiver $Q$, an admissible order $>$ on
$\calB$, and a finite generating set $\mathcal F$ for an ideal $I$ in
$kQ$. Set $R=kQ$ and $\Lambda=kQ/I$. We find a tip-reduced reduced
Gr\"obner basis of uniform elements for $I$ with respect to $>$ and
compute $\rtG$ which must be finite. We also use the sets $O(p)$ and
$N(p)$ defined in Section \ref{section:mainstep}.  We note that, by
the assumption that $\rtG$ is a finite set, therefore $X(p)$ is a
finite set and hence both $O(p)$ and $N(p)$ are finite sets.

We input $M$ in the algorithm as a matrix. In particular, suppose that 
\[\amalg_{i=1}^dw_i\Lambda \extto{\varphi}
\amalg_{i=1}^{t_0}v_i\Lambda \to M\to 0\]
is projective $\L$-presentation of $M$. Then we represent $M$ as the
matrix $(s_{ij})_{i=1,j=1}^{t_0,d}$, where
$\varphi(w_j)=(s_{1j},\ldots,s_{t_0j})$. Note that $s_{ij}$ is in
$v_i\L w_j$.\medskip

\noindent {\bf INPUT}: Nonnegative integers $N$, $t_0$ and $d$,
vertices $v_1,\dots,v_{t_0}$, $w_1,\dots,w_d$, and a $t_0\times
d$-matrix $D$ whose $(i,j)$-th entry is in $v_i\Lambda w_j$.\medskip
 
\noindent {\bf OUTPUT}: For $0\le n\le N$, nonnegative integers $t_n$
and $u_n$, $\{f^i_n\}_{i=1}^{t_n}$, $\{{f^i_n}'\}_{i=1}^{u_n}$ and, if
$n\ge 1$, $h^{n-1,n}_{a,b}$ for $1\le a \le t_{n-1}$ and $1\le b\le
t_n$ as in Section \ref{section:preliminaries}.
\begin{enumerate} 
\item[1.] Set $u_0=0$.  {\bf LiftPresentation}($D$) outputs $t_0$,
$t_1$, $u_1$, $\{f^0_i\}_{i=1}^{t_0}$, $\{f^1_i\}_{i=1}^{t_1}$,
$\{{f^1_i}'\}_{i=1}^{u_1}$.  {\bf
CreateMatrix}($\{f^1_i\}_{i=1}^{t_1}$,$\{f^0_i\}_{i=1}^{t_0}$) outputs
$(h^{0,1}_{i,j})$.
\item[2.] Set $j=1$.
\item[3.] While ($j < N$)
\begin{enumerate}
\item[3.1] Let $T_{j+1}=\{(i,q) \mid  1\le i\le t_{j},\text{ and
}q\in O(\tippath(f^{j}_i)\}$ and $t_{j+1}=|T_{j+1}|$. Choosing
$(i,q)\in T_{j+1}$, one at a time, indexing by $l=1,\dots,t_{j+1}$,
output
\[f^{j+1}_l= f^j_ip-\text{\bf
FirstPart}((f^j_ip-\varepsilon_{i^*}cq'g^2_u),\{f^{j}_i\},\{{f^j_i}'\}),
\]
where $q=\tippath(f^j_i)p=q'\tip(g^2_u)$ and
$c=\frac{\coefftip(f^j_i)}{\coefftip(g^2_u)}$ in $k$. 
\item[3.2] Let $U_{j+1}=\{(i,q) \mid  1\le i\le t_{j},\text{ and
}q\in N(\tippath(f^{j}_i)\}$ and $u_{j+1}=|U_{j+1}|$. Choosing
$(i,q)\in U_{j+1}$, one at a time, indexing by $l=1,\dots,u_{j+1}$,
output
\[
{f^{j+1}_l}'= f^{j}_izg^2_u,
\]
where $q=\tippath(f^j_i)z\tip(g^2_u)$.
\item[3.3] {\bf
  CreateMatrix}($\{f^{j+1}_i\}_{i=1}^{t_{j+1}}$,$\{f^j_i\}_{i=1}^{t_j}$)
outputs $(h^{j,j+1}_{a,b})$.
\item[3.4] $j+1 \leftarrow j$.
\end{enumerate}
\end{enumerate}\medskip
The above algorithm outputs the $f^n_i$'s, the ${f^n_i}'$'s and the
$h^{n,n-1}_{ji}$'s. Next we note that reducing an element $x$ of $R$ by
$\calG$ uses a noncommutative division algorithm \cite{G2}.  The
output of this algorithm is called the \emph{normal form of $x$},
which we denote by $\norm(x)$.  We now obtain the desired first $N$
steps of a projective $\Lambda$-resolution of the cokernel of
$\varphi\colon \amalg_{i=1}^dw_i\Lambda \to
\amalg_{i=1}^{t_0}v_i\Lambda$ for the algorithm above as follows.
Since each $f^i_j$ is right uniform, we let $v^i_j$ be the vertex so
that $f_j^i=f_j^iv^i_j$.  For $n=0,\dots,N$, let
$L^n=\amalg_{i=1}^{t_n}v^n_i\Lambda$ and, for $n=1,\dots,N$, define
$e^n\colon L^n\to L^{n-1}$ by $e^n(v^n_i)$ is $\norm(h^{n-1,n}_{j,i})$
in the $v^{n-1}_j$-th component.

Summarizing, we have the following result. 

\begin{thm}\label{thm:algorithm}
Let $\Lambda=R/I$ where $R=kQ$ for some quiver $Q$.  Assume that
$\calG$ is a tip-reduced Gr\"obner basis for $I$ consisting of uniform
elements and assume further that $\rtG$ is finite.  Let $M$ be a right
$\Lambda$-module such that $M$ has a projective presentation as a
$\Lambda$-module of the form
\[\amalg_{i\in\calI}w_i\Lambda \extto{\varphi}
\amalg_{i\in\calI'}v_i\Lambda\to M\to 0\]
with $\calI$ and $\calI'$ finite sets.

Then there is an algorithm to construct a projective resolution of $M$
over $\L$ associated to the $R$-presentation obtained by our algorithm
\emph{\textbf{LiftPresentation}}. 
\end{thm}

\section{Resolutions of linear modules over Koszul algebras}

In this section we modify our construction to produce minimal
projective resolutions of linear modules over Koszul algebras. We
obtain an algorithm to do this and point out that the assumption that
the Gr\"obner basis is finite is no longer needed. For completeness,
we provide some background.

Recall that if $R=kQ$ and $I$ is an ideal generated by length
homogeneous elements, then the length grading on $R$ induces a positive
$\mathbb{Z}$-grading on $\L=R/I$; namely, $\L=\L_0\amalg \L_1\amalg
\L_2\amalg\cdots$, where $\L_0$ is isomorphic to a finite product of
copies of $k$. Let $\rrad=\amalg_{i\geq 1}\L_i$, which is the graded
radical of $\L$. If $M$ is a graded right $\L$-module with $M_n=0$ for
$n\ll 0$, then $M$ has a minimal graded projective $\L$-resolution
$(\mathbb{L},e)$, where minimal means that $e^n(L^n)\subseteq
L^{n-1}\rrad$ for all $n\geq 1$. We say that $(\mathbb{L},e)$ is a
\emph{linear} resolution and that $M$ is a \emph{linear} module if,
for each $n\geq 0$, the graded module $L^n$ is finitely generated in
degree $n$. The algebra $\L$ is a \emph{Koszul algebra} if $\L_0$ is a
linear right $\L$-module. Koszul algebras were introduced in \cite{P}
and we refer the reader to \cite{BGS,GM,GM2} for further details. Let
$J$ denote the ideal in $kQ$ generated by the arrows.

\begin{thm}
Let $\L=kQ/I$ be a Koszul algebra with $I$ in $J^2$, and let $M$ be a
linear right $\L$-module. Suppose that a start of a minimal projective
linear resolution
\[\amalg_{i=1}^{t_1}w_i\L \to \amalg_{i=1}^{t_0}v_i\L \to M\to 0\]
is given for $M$, where $v_i$ and $w_i$ are vertices in $Q$ and $t_0$
and $t_1$ are positive integers.
\begin{enumerate}
\item[(a)] Then there exists a projective presentation 
\[0\to (\amalg_{i=1}^{t_1} f^1_iR)\amalg (\amalg_{j\in
U_1} {f^1_j}'R)\extto{H^1}\amalg_{i=1}^{t_0} f^0_iR\extto{\pi} M\to
0\] of $M$ as a right $R$-module, where the elements $\{f^0_i\}_{i\in
T_0}$ are vertices, the sets $\{f^1_i\}_{i\in T_1}$ and
$\{{f^1_i}'\}_{i\in U_1}$ are right uniform and right tip-reduced and
contained in $\amalg_{i=1}^{t_0} f^0_iR$, and can be chosen such that
every coordinate of each $f^1_i$ as an element in $\amalg_{i=1}^{t_0}
f^0_iR$ is a sum of elements of length $1$ in $R$ and each ${f^1_i}'$
is in $\amalg_{i\in T_0} f^0_iI$.
\item[(b)]\sloppy There is an algorithm to construct a finite set of elements
$\{f^2_i\}_{i\in T_2}$ in $\amalg_{i=1}^{t_1} f^1_iR$ with
  $f^2_i=\sum_{l=1}^{t_1} f^1_l r_l$ for some linear elements $r_l$ in
  $R$ such that 
\[\amalg_{i\in T_2} f^2_iR/\amalg_{i\in T_2} f^2_iI\extto{e^2}
\amalg_{i=1}^{t_1}f^1_iR/\amalg_{i=1}^{t_1}f^1_iI \to
\Omega_\L^1(M)\to 0\] 
is a start of a minimal projective linear resolution of
$\Omega_\L^1(M)$, where $\Omega^1_\L(M)$ is
$\Ker(\amalg_{i=1}^{t_0}v_i\L \to M)$ and the map $e^2$ is induced by
the inclusion $\amalg_{i\in T_2} f^2_iR\hookrightarrow
\amalg_{i=1}^{t_1}f^1_iR$ as in our earlier construction.
\end{enumerate}
\end{thm}
\begin{proof} (a) The presentation 
\[\amalg_{i=1}^{t_1}w_i\L \to \amalg_{i=1}^{t_0}v_i\L \to M\to 0\]
of $M$ gives rise to the exact sequence
\[0\to K\extto{\varphi} \amalg_{i=1}^{t_0}v_iR \to M\to 0\]
of right $R$-modules. Then $K$ is a projective $R$-module which maps
onto $\Omega^1_\L(M)$. It is easy to see that the natural map
$\amalg_{i=1}^{t_1}w_iR\to \Omega_\L^1(M)$ is a projective cover in
the category of graded right $R$-modules and degree $0$ homomorphisms,
hence there are degree zero maps $\alpha\colon
\amalg_{i=1}^{t_1}w_iR\to K$ and $\beta\colon K\to
\amalg_{i=1}^{t_1}w_iR$ such that
$\beta\alpha=\id_{\amalg_{i=1}^{t_1}w_iR}$. In particular
\[K=\Im\alpha\amalg \Ker\beta \simeq (\amalg_{i=1}^{t_1}w_iR)\amalg
\Ker\beta.\] 

Since $\amalg_{i=1}^{t_0} v_iI$ is the kernel of the map
$\amalg_{i=1}^{t_0} v_iR\to \amalg_{i=1}^{t_0} v_i\L$, we have that
$\Ker\beta$ is contained in $\amalg_{i=1}^{t_0}v_iI$. As $\Ker\beta$
is a projective $R$-module, there are vertices $\{w_i'\}_{i\in U_1'}$
in $Q_0$ for some index set $U_1'$ such that $\Ker\beta\simeq
\amalg_{i\in U_1'} w_i'R$.

Now let 
\begin{xalignat}{2}
f^0_i & = v_i && \text{for $i=1,\ldots,t_0$,}\notag\\
h^1_i & = \varphi\alpha(w_i) && \text{for $i=1,\ldots,t_1$}\notag
\intertext{and}
{h^1_i}' & =\varphi(w_i') && \text{for $i\in U_1'$.}\notag
\end{xalignat}
Since $w_i$ and $w_i'$ are vertices in $Q$, the elements $h^0_1$ and
${h^1_i}'$ are clearly right uniform.  Right tip-reduce each of the
sets $\{h^1_i\}_{i=1}^{t_1}$ and $\{{h^1_i}'\}_{i\in U_1'}$, and
denote the result by $\{f^1_i\}_{i=1}^{t_1}$ and $\{{f^1_i}'\}_{i\in
U_1}$, respectively. The elements are still right uniform.

Since the map $\alpha$ has degree zero and $M$ is a linear
$\L$-module, each of the coordinates of the elements
$\{f^1_i\}_{i=1}^{t_1}$ as elements in $\amalg_{i=1}^{t_0} f^0_iR$ are
all a sum of elements of length $1$ in $R$. The elements
$\{{f^1_i}'\}_{i\in U_1}$ are in $\amalg_{i=1}^{t_0} f^0_iI$, so that
each of the coordinates of an element ${f^1_i}'$ as an element of
$\amalg_{i=1}^{t_0} f^0_iR$ is a sum of elements of length at least
$2$ in $R$. This completes the proof of (a). 

(b) First we look at the construction of $f^2_i$'s given in Section
\ref{section:mainstep}. By linearity, all the coordinates of the
elements $f^2_i$'s occurring in a minimal projective linear resolution
of $M$ as elements in $\amalg_{i=1}^{t_0} f^0_iR$ are a sum of
elements of length $2$ in $R$. Since $I$ is generated by length
homogeneous elements of degree $2$, there is a tip-reduced uniform
Gr\"obner basis consisting of length homogeneous elements of degrees
at least $2$.  An element $g^2_j$ of degree $d$ in $\calG$, occurring
in the construction of a $f^2_i$, gives rise to a homogeneous $f^2_i$
of degree $d$. Tip-reduction does not change the degree, so that to
obtain all the $f^2_i$'s to continue the minimal projective linear
resolution of $M$, we only need to consider the elements of degree $2$
in $\calG$. There are only a finite number of such elements, since
$\calG$ is tip-reduced.

Let $s=(i,q)$ be in $T_2$. Then there
is a $j$ such that
\[q=\tippath(f^1_i)p=q'\tip(g^2_j)\]
for some paths $p$ and $q'$ and $g^2_j$ in $\mathcal{G}$ is the end
relation of $q$, where $\tip(g^2_j)$ and $\tippath(f^1_i)$
overlap. Suppose $\tip(g^2_j)$ is a path of length $2$. It overlaps
$\tippath(f^1_i)$, hence $p$ must be a path of length $1$ (an arrow),
and $q'$ is a vertex. The element $f^1_ip-\varepsilon_{i^*}cq'g^2_j$
is in $(\amalg_{i\in T_1} f^1_iR)\amalg (\amalg_{j\in U_1} {f^1_j}'R)$
with $c=\frac{\coefftip(f^1_i)}{\coefftip(g^2_j)}$ in $k$, but since
the set $\{f^1_i\}_{i\in T_1}\cup \{{f^1_i}'\}_{i\in U_1}$ is not
necessarily right tip-reduced there is no apparent algorithm to
express $f^1_ip-\varepsilon_{i^*}cq'g^2_j$ in this direct sum. Since
all the coordinates of this element have degree $2$ as an element in
$\amalg_{i\in T_0} f^0_iR$, an element ${f^1_i}'$ of degree at least
$3$ does not occur in this expression, so that we only need to
consider the ${f^1_i}'$'s of degree $2$. Since there is a finite
number of paths of length $2$ and since the set of elements of
homogeneous degree $2$ in $\{{f^1_i}'\}_{i\in U_1}$ is right
tip-reduced, there are only a finite number of such elements. Say they
are $\{{f^1_i}'\}_{i\in U_1(2)}$ for some set finite $U_1(2)$. 

In the construction in Section \ref{section:mainstep} of $f^2_s$ we
are assuming that the set $\{f^1_i\}_{i\in T_1}\cup \{{f^1_i}'\}_{i\in
U_1}$ is right tip-reduced. In our case we only have that each of the
sets $\{f^1_i\}_{i=1}^{t_1}$ and $\{{f^1_i}'\}_{i\in U_1(2)}$ is right
tip-reduced, but not necessarily the union. A tip of an ${f^1_j}'$
cannot reduce a tip of an $f^1_i$ by length arguments. So in order to
right tip-reduce the set $\{f^1_i\}_{i\in T_1}\cup \{{f^1_i}'\}_{i\in
U_1}$ the elements $\{f^1_i\}_{i\in T_1}$ stay unchanged, while the
elements $\{{f^1_i}'\}_{i\in U_1}$ might change. We need only the
elements obtained from the set $\{{f^1_i}'\}_{i\in U_1(2)}$. Denote
these new elements by $\{{\mathbf{f}^1_i}'\}_{i\in U_1(2)'}$ for some
finite set $U_1(2)'$, where we record how ${\mathbf{f}^1_i}'$ is
expressed in terms of $f^1_j$'s and ${f^1_j}'$'s. Furthermore the
right tip-reduction of some ${f^1_i}'$ is obtained by subtracting
elements of the form $df^1_ia$, where $a$ is an arrow and $d$ is in
$k$ and elements of the form $d{f^1_j}'$'s, where $j$ is in $U_1(2)$
and $d$ is in $k$. Therefore each ${\mathbf{f}^1_i}'$ is still
homogeneous of degree $2$.

When constructing $f^2_s$ for $s=(i,q)$ in $T_2$ we can
algorithmically find a presentation 
\[f^1_ip-\varepsilon_{i^*}cq'g^2_j=\sum_{l\in T_1} f^1_lr_l+\sum_{l\in
  U_1(2)'} {\mathbf{f}^1_l}'s_l\] 
for some elements $r_l$ and $s_l$ in $R$.  The left hand side has all
coordinates being a sum of elements of length $2$. If some path of
length at least $2$ occurs in some $r_l$, then $\tippath(\sum_{l\in
T_1} f^1_lr_l)$ is equal to $\tippath(\sum_{l\in U_1(2)'}
{\mathbf{f}^1_l}'s_l)$. As we have seen before this contradicts the
fact that $\{f^1_s\}_{s\in T_1}\cup \{{\mathbf{f}^1_{s}}'\}_{s\in
U_1(2)'}$ is right tip-reduced. By length arguments no vertex can
occur in any $r_l$. Hence each $r_l$ is a sum of elements of length
$1$ and each $s_l$ is a vertex. Substituting ${\mathbf{f}^1_{s}}'$
with the expressions in $f^1_j$'s and ${f^1_j}'$'s, we obtain as
before
\[f^2_s=f^1_ip-\sum_{l\in T_1}
f^1_lr'_l=\varepsilon_{i^*}cq'g^2_j+\sum_{l\in U_1(2)} {f^1_l}'s_l'\]
for some linear elements $r'_l$ in $R$ and some elements $s_l'$ in $R$
of degree $0$, and therefore all the coordinates of $f^2_s$ as
elements in $\amalg_{i=1}^{t_0} f^0_iR$ are a sum of elements of
length $2$.

Now let $x$ be a homogeneous element of degree $2$ in
$(\amalg_{i=1}^{t_1} f^1_iR)\cap (\amalg_{i=1}^{t_0}f^0_iI)$, where
$\tip(x)$ is smallest possible such that $x$ is not in $\amalg_{i\in
T_2} f^2_iR$. Hence $x=\sum_{i=1}^{t_1}f^1_ib_i=\sum_{i=1}^{t_0}
f^0_ig^2_{j_i}$, where $b_i$ is a sum of elements of length $1$ in $R$
and $g^2_{j_i}$ is homogeneous elements of degree $2$ in $\calG$. Then
$\tip(x)=\tip(f^1_ib_i)=f^0_l\tip(g^2_{j_l})$ for some $i$ and some
$l$. Then $\tip(x)$ is equal to $\tip(f^2_j)$ for some $j$ in $T_2$.
By the choice of $x$ the element $x-cf^2_j$ is in $\amalg_{i\in T_2}
f^2_iR$, for some element $c$ in $k$. Hence all homogeneous elements
of degree $2$ in $(\amalg_{i=1}^{t_1} f^1_iR)\amalg
(\amalg_{i=1}^{t_0}f^0_iI)$ are in $\amalg_{i\in T_2}
f^2_iR$. Therefore we have constructed all $f^2_i$'s of degree $2$,
and by construction the elements $\{f^2_i\}_{i\in T_2}$ are right
uniform and right tip-reduced. These elements give rise to the minimal
projective cover of $\Omega_\L^2(M)$. Then there is a natural map
$\amalg_{i\in T_2} f^2_iR/\amalg_{i\in T_2} f^2_iI\to \Omega^2_\L(M)$,
which is a projective cover. Then we have a start of minimal
projective linear resolution of $\Omega_\L^1(M)$
\[\amalg_{i\in T_2} f^2_iR/\amalg_{i\in T_2} f^2_iI\extto{e^2}
\amalg_{i=1}^{t_1}w_i\L\to \Omega_\L^1(M)\to 0,\]
since the elements $p$ and $r_l'$ are linear elements. The
construction of the map $e^2$ is the same as given in Section
\ref{section:mainstep}. 
\end{proof}
To see that the above arguments actually give rise to an algorithm, we
first note that one need not right tip-reduce the whole set
$\{{h^1_i}'\}_{i\in U_1}$, which maybe infinite. We only need to right
tip-reduce those elements of homogeneous degree $2$ in this set. This
subset of $\{{h^1_i}'\}_{i\in U_1}$ can be chosen to be finite, since
the subspace of elements of homogeneous degree $2$ of each $v_iI$ has
basis the elements of homogeneous degree $2$ of $v_i\rtG$, which is a
finite set. We are also not assuming that the Gr\"obner basis $\calG$
is finite. But we only need the homogeneous elements of degree $2$ in
$\calG$, which is a finite set and may be computed by right
tip-reducing a set of right uniform generators of $I$. It follows that
the construction of the elements $f^2_s$ of homogeneous degree $2$ is
algorithmic. 

We remark that there is another method for constructing the $f^2_i$'s
of homogeneous degree $2$. Namely, let $A$ be the $k$-span of $\{
f^1_ia\mid i=1,\ldots, t_1, a\in Q_1\}$ and $B$ be the $k$-span of
$\{f^0_ig\mid g\in \calG \text{\ and length of\ } g = 2\}$. Then one
may use linear algebra to find a basis $\{ b_1,\ldots,b_{t_2}\}$ of
$A\cap B$. Viewing the $b_i$'s as elements in $\amalg_{i=1}^{t_0}
f^0_iR$, the set $\{f^2_i\}_{i=1}^{t_2}$ can be found by right
tip-reducing $\{b_i\}_{i=1}^{t_2}$. This may be a faster way of
finding the $f^2_i$'s than the method presented in the above proof.

In general the construction in Section \ref{section:construction} does
not produce a minimal projective resolution of a linear module.  The
algorithm of this section differs from the algorithm in Section
\ref{section:mainstep} in that one only considers elements of a
Gr\"obner basis of length $2$ when constructing $T_2$.  We illustrate
this in the following example.

\begin{example}
We continue with Example \ref{example:basicIII}. We note that $\L$ is
a Koszul algebra and that $M=v_1\L/\rrad$ is a linear module. We saw
that the resolution constructed by the algorithm in Section
\ref{section:construction} gave a non-minimal projective resolution of
$M$ for the ordering $>_1$. For this ordering, recall that the
Gr\"obner basis $\calG=\{ab-cd,be,cde\}$ for $I$. Referring back to
Example \ref{example:basicII} in Section \ref{section:mainstep} we now
construct $T_2$ using the algorithm described above, that is; only
using $ab-cd$ and $be$ from $\calG$. In this way we only produce
$f^2_1$ as in Example \ref{example:basicII}. We obtain the same
resolution as given in Example \ref{example:basicIII} for the ordering
$>_2$ in this way, and hence producing a minimal projective
resolution of $M$ over $\L$. 
\end{example}

\end{document}